\documentclass[letterpaper, 10 pt, conference]{ieeeconf}  % Comment this line out if you need a4paper

\IEEEoverridecommandlockouts                              % This command is only needed if 
\usepackage{graphics} % for pdf, bitmapped graphics files
\usepackage{epsfig} % for postscript graphics files
\usepackage{mathptmx} % assumes new font selection scheme installed
\usepackage{times} % assumes new font selection scheme installed
\usepackage{amsmath} % assumes amsmath package installed
\usepackage{amssymb}  % assumes amsmath package installed

\usepackage{amssymb} % for mathbb
\usepackage{dsfont} % for mathds
\usepackage{amsmath} % for align
\newtheorem{theorem}{Theorem}[section]
\newtheorem{remark}[theorem]{Remark}
\usepackage{xfrac}
\usepackage{cases} % for subnumcases
\usepackage{xcolor} % for textcolor
\usepackage{cite}

\allowdisplaybreaks

\title{\LARGE \bf
Optimisation of Region of Attraction Estimates for the Exponential Stabilisation of the Intrinsic Geometrically Exact Beam Model}

\author{Marc Artola$^{1}$, Charlotte Rodriguez$^{2}$, Andrew Wynn$^{1}$, Rafael Palacios$^{1}$, and G\"{u}nter Leugering$^{2}$% <-this % stops a space
\thanks{*Marc Artola and Charlotte Rodriguez are members of the Innovative Training Network
ConFlex. This project has received funding from the European
Union's Horizon 2020 research and innovation programme under the Marie Sklodowska-Curie grant agreement No 765579.}% <-this % stops a space
\thanks{$^{1}$Marc Artola, Andrew Wynn and Rafael Palacios are with the Department of Aeronautics, Imperial College London, Exhibition Road, London SW7 2AZ, UK}
\thanks{
        {\tt\small marc.artola16@imperial.ac.uk}}
\thanks{
        {\tt\small a.wynn@imperial.ac.uk}}
\thanks{
        {\tt\small r.palacios@imperial.ac.uk}}%
\thanks{$^{2}$Charlotte Rodriguez and G\"{u}nter Leugering are with the Department of Mathematics Chair of Applied Mathematics, Friedrich-Alexander Universit\"{a}t, Cauerstr. 11 91058 Erlangen}
\thanks{
        {\tt\small charlotte.rodriguez@fau.de}}
\thanks{
        {\tt\small guenter.leugering@fau.de}}%}
    }

\usepackage{tikz}
\newcommand\copyrighttext{%
  \footnotesize \textcopyright 2021 IEEE.  Personal use of this material is permitted.  Permission from IEEE must be obtained for all other uses, in any current or future media, including reprinting/republishing this material for advertising or promotional purposes, creating new collective works, for resale or redistribution to servers or lists, or reuse of any copyrighted component of this work in other works.
}
\newcommand\copyrightnotice{%
\begin{tikzpicture}[remember picture,overlay]
\node[anchor=south,yshift=10pt] at (current page.south) {\fbox{\parbox{\dimexpr\textwidth-\fboxsep-\fboxrule\relax}{\copyrighttext}}};
\end{tikzpicture}%
}

\begin{document}

\maketitle
\copyrightnotice
\thispagestyle{empty}
\pagestyle{empty}

%%%%%%%%%%%%%%%%%%%%%%%%%%%%%%%%%%%%%%%%%%%%%%%%%%%%%%%%%%%%%%%%%%%%%%%%%%%%%%%%
\begin{abstract}
A systematic approach to maximise estimates on the region of attraction in the exponential stabilisation of geometrically exact (nonlinear) beam models via boundary feedback is presented. Starting from recently established stability results based on Lyapunov arguments, the main contribution of the presented work is to maximise the analytically found bounds on the initial datum, for which local exponential stability is guaranteed, via search of (optimal) polynomial Lyapunov functionals using an iterative semi-definite programming approach.
\end{abstract}

%%%%%%%%%%%%%%%%%%%%%%%%%%%%%%%%%%%%%%%%%%%%%%%%%%%%%%%%%%%%%%%%%%%%%%%%%%%%%%%%
\section{INTRODUCTION}

There is a growing interest in beam models describing the three-dimensional motions of highly-flexible light-weight structures -- for instance, robotic arms \cite{grazioso2018robot}, flexible aircraft wings \cite{palacios10aircraft, Artola2020aero} or wind turbine blades \cite{wang2014windturbine} --, which exhibit motions of large magnitude, not negligible in comparison to the overall dimensions of the object. To capture such a behaviour, one needs so-called \emph{geometrically exact} beam models, which then exhibit nonlinearities. 
Such models, similar to the canonical Euler-Bernoulli and Timoshenko models, are one-dimensional with respect to the spatial variable, and account for small strains. The independent variable $x\in [0, \ell]$ ranges along the centerline of the beam, $\ell>0$ being the beam's total arclength.

In engineering applications, there is a clear need to eliminate vibrations and flutter in these structures \cite{Matsuoka1995, UchiyamaKonno1991}. This translates into the task of finding appropriate controls (here we consider boundary feedback control) to make the mathematical model describing the beam exponentially stable, relying on Lyapunov-based arguments and in a sense made clear in the following section. We will see that exponential stability may be achieved at least locally (i.e., for small initial data) due to the nonlinear nature of the system. Therefore, it becomes interesting to determine when (i.e., for which initial states of the system) one shall expect exponential decay of the solutions in the case of a freely vibrating beam -- meaning that external forces such as gravity or aerodynamic forces are set to zero. In order to do so, we turn to semi-definite programming, a numerical approach which is gaining popularity in PDE analysis and control~\cite{GOULART2012692, Valmorbida2016, Marx2020}. This is employed here to systematise the choice of Lyapunov functional, leading to sharper bounds on the recently derived region of attraction estimates for the boundary feedback stabilisation of geometrically exact beams~\cite{RL2020}. 

The paper is organised as follows. Section \ref{sec:formulation} introduces the mathematical formulation for geometrically-nonlinear beams, while in \S~\ref{sec:boundary_feedback} a concise description of the proposed boundary feedback control strategy is given, summarising the previously found stability results. 
%\textcolor{red}{In \S~\ref{sec:optimisation}, the optimisation problem, suitable for semi-definite programming, to maximise the region of attraction estimates is proposed.} 
The main contribution of this work is introduced in \S\ref{sec:optimisation}, where an optimisation problem, solved iteratively using semi-definite programming and designed to maximise the region of attraction estimates, is proposed. Finally, the methodology is applied on a numerical example in \S\ref{sec:results}.
\section*{NOTATION}
Let us introduce some useful notation. We denote by $|\cdot|$ and $\langle \cdot , \cdot \rangle$ the Euclidean norm and inner product in $\mathbb{R}^n$, and for any matrix $M$, $\|M\|$ is the operator norm induced by $|\cdot|$. The identity and null matrices are denoted by $\mathbf{I}_n \in \mathbb{R}^{n \times n}$ and $\mathbf{0}_{n, m} \in \mathbb{R}^{n \times m}$, and we use the abbreviation $\mathbf{0}_{n} = \mathbf{0}_{n, n}$. The transpose of any matrix $M$ is denoted $M^\top$.
For any $M \in \mathbb{R}^{n \times n}$, we say that $M$ is positive (semi-)definite, and denote it $M \succ 0$ (resp. $M \succcurlyeq 0$), if $\langle u, Mu \rangle >0$ (resp. $\geq 0$) for all $u \in \mathbb{R}^n \setminus \{0\}$. 
%For any diagonal matrix $M \in \mathbb{R}^{12}$, we denote $M = \mathrm{diag}(M_-, M_+)$ with $M_-$ and $M_+$ containing the first and last six diagonal entries of $M$, respectively.
The set of (positive definite) diagonal matrices of size $n$ is denoted $\mathbb{D}^n$ (resp. $\mathbb{D}_{++}^n$).
Also, for any $Q(\cdot)$ with values in $\mathbb{D}^{12}$, we denote by $\{q_i(x)\}_{i=1}^{12}$ the diagonal entries of $Q(x)$, and write $Q = \mathrm{diag}(Q_-, Q_+)$ where
\begin{align} \label{eq:Q_diag_entries}
    Q_- = \mathrm{diag}(q_1, \ldots, q_6), \quad Q_+ = \mathrm{diag}(q_7, \ldots, q_{12}).
\end{align}

$\mathrm{SO}(3)$ is the set of unitary real matrices of size $3$ and with a determinant equal to $1$, also called \emph{rotation} matrices.The cross product between any $u, \zeta \in \mathbb{R}^3$ is denoted $u \times \zeta$, and we shall also write $\widehat{u} \,\zeta = u \times \zeta$, meaning that $\widehat{u}$ is the skew-symmetric matrix
\begin{align*}
\widehat{u} = \begin{bmatrix}
0 & -u_3 & u_2 \\
u_3 & 0 & -u_1 \\
-u_2 & u_1 & 0
\end{bmatrix}, 
\end{align*}
while $u$ is recovered by means of the operator $\mathrm{vec}(\cdot)$ acting on skew-symmetric matrices as follows: $\mathrm{vec}(\widehat{u}) = u$.
% and for any skew-symmetric $\mathbf{u} \in \mathbb{R}^{3 \times 3}$, the vector $\mathrm{vec}(\mathbf{u}) \in \mathbb{R}^3$ is such that $\mathbf{u} = \widehat{\mathrm{vec}(\mathbf{u})}$.

\section{PROBLEM FORMULATION}
\label{sec:formulation}
Commonly, the mathematical model for geometrically exact beams is a quasilinear second-order system written in terms of the position of the beam's centerline, $\mathbf{p}(x,t) \in \mathbb{R}^3$, and the orientation of its cross sections given by the columns of the matrix $\mathbf{R}(x,t) \in \mathrm{SO}(3)$, both being expressed in some fixed coordinate system such as the standard basis $\{e_i\}_{i=1}^3$ of $\mathbb{R}^3$. It is set in $(0, \ell)\times(0, T)$ and reads
\begingroup % keep the change local
\setlength\arraycolsep{2pt}
\begin{align*}
% \partial_t \left[ \begin{bmatrix}
% \mathbf{R} & \mathbf{0}_{3}\\ \mathbf{0}_{3} & \mathbf{R}
% \end{bmatrix} \mathbf{M}
% v
% \right] = \partial_x \left[ \begin{bmatrix}
% \mathbf{R} & \mathbf{0}_{3}\\ \mathbf{0}_{3} & \mathbf{R}
% \end{bmatrix} \mathbf{C}^{-1} s\right] + 
% \begin{bmatrix}
% \mathbf{0}_3  & \mathbf{0}_3 \\ (\partial_x \widehat{\mathbf{p}})\mathbf{R} & \mathbf{0}_3
% \end{bmatrix}
% \mathbf{C}^{-1}s,
\begin{bmatrix}
\partial_t & \mathbf{0}\\
(\partial_t \widehat{\mathbf{p}}) & \partial_t
\end{bmatrix} \left[ \begin{bmatrix}
\mathbf{R} & \mathbf{0}\\ \mathbf{0} & \mathbf{R}
\end{bmatrix}
\mathbf{M} v \right] = \begin{bmatrix}
\partial_x & \mathbf{0} \\ (\partial_x \widehat{\mathbf{p}}) & \partial_x
\end{bmatrix} \left[\begin{bmatrix}
\mathbf{R} & \mathbf{0}\\ \mathbf{0} & \mathbf{R}
\end{bmatrix} \mathbf{C}^{-1} s \right],
\end{align*}
\endgroup
where $\mathbf{M},\mathbf{C} \in \mathbb{R}^{6\times 6}$ are the so-called \emph{mass} and \emph{ flexibility} matrices, $v(x,t) \in \mathbb{R}^6$ contains the linear and angular \emph{velocities}, and $s(x,t) \in \mathbb{R}^6$ contains the linear and angular \emph{strains}:
\begin{equation*} %\label{eq:v_s}
v
= \begin{bmatrix}
\mathbf{R}^\top \partial_t \mathbf{p}\\ \mathrm{vec}\left( \mathbf{R}^\top \partial_t \mathbf{R} \right)
\end{bmatrix}, \quad 
s= \begin{bmatrix}
\mathbf{R}^\top \partial_x \mathbf{p}  - e_1 \\ 
\mathrm{vec}\left(\mathbf{R}^\top \partial_x \mathbf{R} \right) - \Upsilon_c
\end{bmatrix},
\end{equation*}
in which $\Upsilon_c = \mathrm{vec}( R^\top \frac{\mathrm{d}}{\mathrm{d}x}R )$ is the curvature before deformation, the given matrix $R(x) \in \mathrm{SO}(3)$ describing the cross sections' orientation before deformation. 
This system is the Geometrically Exact Beam model, due to Reissner \cite{reissner1981finite}, who initially derived the static formulation, and Simo \cite{simo1985finite}, who extended it to the dynamic case.

The mathematical model may also be written in terms of so-called \emph{intrinsic} variables: the velocities $v$ and strains $s$ (or equivalently velocities and internal forces and moments), which are all expressed in a moving basis attached to the beam's centerline -- namely, the basis defined by the columns of $\mathbf{R}$. This yields the Intrinsic Geometrically Exact Beam model, 
%\textcolor{red}{(Should we rather call it the ``intrinsic beam model'' as in you paper? For me it does not matter.)}\textcolor{blue}{Fine for me to keep with your nomenclature}, 
or IGEB, due to Hodges \cite{hodges2003geometrically}, which reads
\begin{align} \label{eq:IGEB_pres}
\partial_t y + A \partial_x y + \overline{B}(x)y = \overline{g}(y),
\end{align}
the unknown being $y = [v^\top, s^\top]^\top$. Boundary conditions for \eqref{eq:IGEB_pres} can be generally expressed as $G_\partial y(x_\partial)=b_\partial$, with $G_\partial\in\mathbb{R}^{6\times12}$, $b_\partial\in\mathbb{R}^{6}$ and $x_\partial=\{0,\ell\}$ denoting each of the two boundaries of a beam of arclength $\ell$. We shall henceforth focus solely on this formulation.
An advantageous feature of this model is that it falls into the class of \emph{one-dimensional first-order hyperbolic systems} (hyperbolic meaning that $A$ has real eigenvalues only and twelve associated independent eigenvectors) and thus provides access to a broad mathematical literature beyond the context of beam models (e.g., \cite{Li_Duke85, BC2016}).
Moreover, the IGEB model is only \emph{semilinear}, with the nonlinear function $\overline{g}$ being quadratic, thus locally Lipschitz. It is defined by
\begin{align}\label{eq:gbar}
\overline{g}(y) = -
\begin{bmatrix}
\mathbf{M}^{-1}L_1(v) & \mathbf{M}^{-1}L_2(\mathbf{C}^{-1}s) \\
\mathbf{0}_6 & -L_1(v)^\top
\end{bmatrix}
\begin{bmatrix}
\mathbf{M}v \\ s
\end{bmatrix},
\end{align}
where for any $u = [u_1^\top, u_2^\top]^\top$ with $u_i \in \mathbb{R}^3$, 
\begin{align*}
L_1(u) = \begin{bmatrix}
\widehat{u}_2 & \mathbf{0}_3 \\
\widehat{u}_1 & \widehat{u}_2
\end{bmatrix}, \quad L_2(u) = \begin{bmatrix}
\mathbf{0}_3 & \widehat{u}_1\\
\widehat{u}_1 & \widehat{u}_2
\end{bmatrix}.
\end{align*}
The matrix $\overline{B}(x)$ may not be assumed arbitrarily small, hence, not only is the linearised system not homogeneous, but also \eqref{eq:IGEB_pres} cannot be seen as the perturbation of a system of conservation laws, which makes the stabilisation study challenging. Denoting $\mathbf{E} = L_1(\overline{s})$ with $\overline{s} = [e_1^\top, \Upsilon_c^\top]^\top$, one has
\begin{align*}
A = \begin{bmatrix}
\mathbf{0}_6 & -(\mathbf{M}\mathbf{C})^{-1} \\
-\mathbf{I}_6 & \mathbf{0}_6
\end{bmatrix}, \quad \overline{B} = 
\begin{bmatrix}
\mathbf{0}_{6} & -\mathbf{M}^{-1}\mathbf{E}\mathbf{C}^{-1} \\
\mathbf{E}^\top & \mathbf{0}_{6} 
\end{bmatrix}.
\end{align*}
%
%\textcolor{red}{(Don't hesitate to remove/point out references that you think I should remove)}
The neater structure of the IGEB model makes it well-suited for aeroelastic modelling and engineering, notably in the context unmanned aircraft aiming to remain airborne over large time horizons, which are consequently very light-weight and slender and exhibit great flexibility \cite{palacios10aircraft, Palacios2011intrinsic}. 
% (see also \cite{Artola2020aero, Artola2021damping} where the authors additionally take into account structural damping).
Furthermore, one may also see the IGEB model as the beam dynamics formulated in the \emph{Hamiltonian} framework (see \cite[Sections 5, 6]{Simo1988}), leading to the study of this system from the perspective of \emph{Port-Hamiltonian Systems} when taking into account the interactions of the beam with its environment (see \cite{Macchelli2009}).
%and \cite[Section 4.3.2]{Macchelli2009book}).

\begin{remark}
We suppose that $\mathbf{M}$ and $\mathbf{C}$ are independent of $x$, and restrict our study to beams made of an isotropic material, with constant parameters, and such that sectional principal axes are aligned with the body-attached basis:
\begin{align*}
\mathbf{M} &= \rho \mathrm{diag}( a,  a,  a, (I_2+I_3)k_1, I_2, I_3),\\
\mathbf{C} &= \rho \mathbf{M}^{-1}\mathrm{diag}(E, k_2G, k_3G, G, E, E)^{-1},
\end{align*}
with density $\rho>0$, cross section area $a>0$, shear modulus $G>0$, Young modulus $E>0$, area moments of inertia $I_2, I_3>0$, shear correction factors $k_2, k_3 >0$, and the factor $k_1>0$ that corrects the polar moment of area.
\end{remark}

\section{BOUNDARY FEEDBACK STABILISATION}
\label{sec:boundary_feedback}

To the best of our knowledge, global in time existence and uniqueness of $C^0([0, +\infty), H^1(0, \ell; \mathbb{R}^{12}))$ (or $C^0([0, \ell  ] \times  [0, +\infty ); \mathbb{R}^{12})$) solutions to \eqref{eq:IGEB_pres} is not provided by general results present in the literature.
However, for some specific closed-loop problems for \eqref{eq:IGEB_pres}, one may deduce well-posedness on the entire time interval, together with an exponential decay in time of the solution. 
Being then based on local in time solutions to \eqref{eq:IGEB_pres} and on maintaining the nonlinear term to be small throughout the proof, this stability result is \emph{local}, in the sense that it holds for small initial data only.

More precisely, we can stabilise the beam by applying a \emph{velocity feedback control} of the form $-\overline{\kappa}v$ at the end $x=0$ (i.e., effectively a damper, where the force at the boundary is constrained to oppose velocity), with $\overline{\kappa} \in \mathbb{R}^{6\times 6}$, while the other end is clamped. This yields the system
%\begin{subequations}
\begin{align}\label{eq:IGEB}
\left\{ 
\begin{aligned}
&\partial_t y + A \partial_x y + \overline{B}(x)y = \overline{g}(y) &&\text{in } (0,\ell)\times(0,T) \\
&v(\ell, t) = 0 &&\text{for } t \in (0,T)\\ 
&-\mathbf{C}^{-1} s(0, t) = -\overline{\kappa} v(0, t) && \text{for } t \in (0,T)\\
&y(x,0) = y^0(x)  && \text{for } x \in (0,\ell),
\end{aligned}
\right.
\end{align}
with initial datum $y^0 = ({v^0}^\top, {s^0}^\top)^\top$ in $\mathbf{H}^1$, where we henceforth use the shortened notation $\mathbf{H}^1 = H^1(0, \ell; \mathbb{R}^{12})$. Then, one has the following theorem.

\begin{theorem}
\label{thm:stabilization}
Assume that $R \in C^1([0, \ell]; \mathrm{SO}(3))$ and $\overline{\kappa} \in \mathbb{R}^{6\times 6}$ is symmetric positive definite. Then, the steady state $y \equiv 0$ of \eqref{eq:IGEB} is \emph{locally $H^1$ exponentially stable}, in the sense that there exist $\varepsilon>0$, $\alpha>0$ and $\eta \geq 1$ such that, for all $y^0 \in \mathbf{H}^1$ satisfying $\|y^0\|_{\mathbf{H}^1}\leq \varepsilon$ and the compatibility conditions $v^0(\ell) = 0$, $\mathbf{C}^{-1} s^0(0) = \overline{\kappa} v^0(0)$, there exists a unique global in time solution $y \in C^0([0, +\infty); \mathbf{H}^1)$ to \eqref{eq:IGEB}, and
%$\|y(\cdot, t)\|_{\mathbf{H}^1} \leq \eta  e^{- \alpha t } \|y^0\|_{\mathbf{H}^1}$ for all $t \in [0, +\infty)$.
\begin{align*}
\|y(\cdot, t)\|_{\mathbf{H}^1} \leq \eta  e^{- \alpha t } \|y^0\|_{\mathbf{H}^1}, \quad \text{for all }t \in [0, +\infty).
\end{align*}
\end{theorem}

\medskip

Theorem \ref{thm:stabilization} is proved in \cite[Th. 1.5]{RL2020} with more constraints on $\overline{\kappa}$, but the proof is in fact easily adjusted to any positive definite symmetric $\overline{\kappa}$; see also \cite[Th. 2.4, Rem. 2]{R2021} for the case of a general linear elastic material law (anisotropic material, beam parameters dependent on $x$).
Let us now give some detail on the idea of the proof.
% -- where IGEB model is written in terms of velocities and internal forces and moments.

\subsection{Riemann Invariants} 
The fact that $A$ is hyperbolic allows us to apply a change of variable and thereby write \eqref{eq:IGEB} in terms of the so-called \emph{Riemann invariants} (or \emph{diagonal} or \emph{characteristic} form), a form in which the stability analysis is simplified.
Applying the change of variables
\begin{align*}
r = L y \quad \text{with } \ L = \left[\begin{smallmatrix}
\mathbf{I}_6 & D \\
\mathbf{I}_6 & -D
\end{smallmatrix}\right], \ \ D = (\mathbf{M}\mathbf{C})^{-\sfrac{1}{2}},
\end{align*}
to \eqref{eq:IGEB}, for $\mathbf{D} = \mathrm{diag}(-D, D)$, $B = L \overline{B}L^{-1}$, $g(r) = L \overline{g}(L^{-1}r)$, and $\kappa = (\mathbf{M}D + \overline{\kappa})^{-1} (\mathbf{M}D - \overline{\kappa})$, one obtains
\begin{align}\label{eq:IGEBD}
\left\{ 
\begin{aligned}
&\partial_t r + \mathbf{D} \partial_x r + B(x)r = g(r) &&\text{in } (0,\ell)\times(0,T) \\
&r_-(\ell, t) = - r_+(\ell,t) &&\text{for } t \in (0,T)\\ 
&r_+(0, t) = \kappa r_-(0, t) && \text{for } t \in (0,T)\\
&r(x,0) = r^0(x)  && \text{for } x \in (0,\ell),
\end{aligned}
\right.
\end{align} 
with initial datum $r^0 = Ly^0$. In line with the sign of the diagonal entries of $\mathbf{D}$, we denote $r = (r_-^\top, r_+^\top)^\top$ with $r_-$ and $r_+$ being the first and last six components of $r$, respectively.

From the relationship between $g$ and $\overline{g}$, and the definition \eqref{eq:gbar} of the latter, one can deduce that the components of $g$ also write as $g_i(r) = \langle r, G^i r \rangle$ where each $G^i \in \mathbb{R}^{12 \times 12}$ is a specific symmetric matrix dependent on the beam parameters. We then define the constant $C_g>0$ by $C_g = \sqrt{\sum_{i=1}^{12}\|G^i\|^2}$. We also define $C_B>0$, characterising the linear lower order term, by $C_B = \max_{x \in [0, \ell]}\|B(x)\|$.

\subsection{Lyapunov Functional} \label{subsec:lyap}
One may equivalently study this stability problem for \eqref{eq:IGEB} or for its diagonal form \eqref{eq:IGEBD}.
Among the methods commonly used to study stability, a so-called $H^1$ \emph{quadratic Lyapunov functional} is used in \cite{RL2020, R2021}, namely a functional of the form 
\begin{align}
\mathcal{L}(t) = \sum_{j=0}^1 \int_0^\ell \left\langle \partial_t^j r(x,t) \,, Q(x) \partial_t^j r(x,t) \right\rangle dx
\label{eq:Lyapunov_functional}
\end{align}
for some $Q=Q(x) \in \mathbb{R}^{12 \times 12}$ and for $t \in [0, T]$, which is equivalent to the squared $\mathbf{H}^1$ norm of $r(\cdot, t)$ and has an exponential decay with respect to time, as long as the solution $r$ to \eqref{eq:IGEBD} is in some ball of $C^{0}([0,\ell ]\times[0,T];\mathbb{R}^{12})$. 

For one-dimensional first-order hyperbolic systems, Bastin and Coron \cite{bastin2017exponential} have systematised the search of such functionals and given  sufficient criteria for their existence. More precisely, it is sufficient to find a matrix-valued function $Q$ fulfilling a series of matrix inequalities that involve both the coefficients appearing in the governing equations and the boundary conditions (and thus, the feedback control). For System \eqref{eq:IGEB}, their result yields the following theorem (see \cite[Prop. 3.1]{RL2020}).

\begin{theorem} \label{th:bastin_coron}
If there exists $Q \in C^1([0,\ell ]; \mathbb{D}_{++}^{12})$ with
\begin{subequations}\label{eq:MI}
\begin{align}
\label{eq:MI_ell}
Q_+(\ell ) - Q_-(\ell ) \succcurlyeq 0,\\
\label{eq:MI_0}
Q_-(0)D - \kappa^\top Q_+(0)D \kappa \succcurlyeq 0, \\ 
\label{eq:MI_S}
-\tfrac{\mathrm{d}}{\mathrm{d}x}Q(x) \mathbf{D} + Q(x)B(x) + B(x)^\top Q(x) \succ 0,
\end{align}
\end{subequations}
for any $x \in [0, \ell ]$, then the steady state $y \equiv 0$ of \eqref{eq:IGEB} is locally $H^1$ exponentially stable.
\end{theorem}

For any given $Q$ as in Theorem \ref{th:bastin_coron}, we define the constants $C_\mathcal{S}, C_Q>0$ as follows:
\begin{align*}
C_Q = \max \left\{\max_{ {x \in [0,\ell ], \, { 1\leq i \leq 12 }}} q_i(x), \Big( \min_{ {x \in [0,\ell ],\, {1\leq i \leq 12}}} q_i(x)\Big)^{-1} \right\}
\end{align*}
and, for $\mathcal{S}:=-\tfrac{\mathrm{d}}{\mathrm{d}x}Q \mathbf{D} + QB + B^\top Q$,
\begin{align*}
C_\mathcal{S} = \min_{x \in [0, \ell]} \sigma_{\min}(\mathcal
{S}(x))
\end{align*}
where $\sigma_{\min}(\mathcal{S}(x))$ denotes the smallest eigenvalue of $\mathcal{S}(x)$. 
%Moreover, we introduce the constant $C_Q>0$ defined by 
%\begin{align*}
%C_Q = \max \left\{\max_{ {x \in [0,\ell ], \, { i \in \{1, \ldots, 12\}}}} q_i(x), \Big( \min_{ {x \in [0,\ell ],\, {i \in \{1, \ldots, 12\}}}} q_i(x)\Big)^{-1} \right\},
%\end{align*}
%where $q_i$ denotes the $i$-th diagonal entry of $Q$.

\subsection{Region of Attraction} Going through the proof of Theorem \ref{th:bastin_coron} while keeping track of the constants, one can show that for any given $Q$ fulfilling the assumptions of Theorem \ref{th:bastin_coron}, the constants $\varepsilon,\alpha>0$ and $\eta\geq 1$ introduced in Theorem \ref{thm:stabilization} may be chosen as follows.
The bound on the initial datum has the form
\begin{align} 
\varepsilon = \|L\|^{-1} \min \left\{ \frac{\delta}{2 C_1 \eta}, \frac{\delta_0}{\eta}\right\},
\label{eq:analytical_bound}
\end{align}
where $C_1>0$ is the constant coming from the standard Sobolev inequality $\|\varphi\|_{C^0([0, \ell ]; \mathbb{R}^{12})} \leq C_1 \|\varphi\|_{\mathbf{H}^1}$ (see \cite[Th. 5, Sec. 5.6]{evansPDE}). 
%
%\textcolor{blue}{(Do something about the constant $\delta_0$: need to check that $\delta = \delta_0$ is ok, and/or explain what $\delta_0$ is.)}
In fact, the proof of Theorem \ref{th:bastin_coron} consists in first showing that, at least on a small time interval, there exists a unique solution to \eqref{eq:IGEBD} (see \cite[Th. 10.1]{bastin2017exponential}), before extending this solution for all times by means of the Lyapunov functional $\mathcal{L}$ (see \eqref{eq:Lyapunov_functional}). While the appearance of the constant $\delta_0>0$ is due to the aforementioned existence and uniqueness result, the constant $\delta>0$ is directly related to the decay of $\mathcal{L}$ and, thereby, to the feedback control. Therefore, the latter will be the main focus of our work. To be more precise, $\delta$ can be chosen as any \emph{positive} number such that
\begin{align} \label{eq:alpha}
\alpha := \tfrac{1}{2}C_Q\big(C_\mathcal{S} - 4 C_Q C_g  \delta \big) >0
\end{align}
holds, where $\alpha>0$ is then the exponential decay of the solution. 
Note that there is a competition in \eqref{eq:alpha} between the exponential decay and $\delta>0$, and that the constant $C_\mathcal{S}$ depends on $Q$.
Finally, $\eta\geq 1$ is given by
\begin{align} \label{eq:eta}
\eta = C_Q\big(2 C_\mathrm{IB}(\delta) + 1 \big) \|L\| \|L^{-1}\|
\end{align}
where $C_\mathrm{IB}(\delta)>0$ is defined by
\begin{align} \label{eq:def_CIB}
C_\mathrm{IB}(\delta) = \max \left\{ \|\mathbf{D}\|, \|\mathbf{D}\|^{-1}, (C_B + C_g \delta), \frac{C_B + C_g \delta}{\|\mathbf{D}\|} \right\}.
\end{align}

In view of this, our aim in what follows, is to go beyond the specific Lyapunov functional found in \cite{RL2020, R2021}, and look for an optimal functional, in the form of polynomials, such that the bound constraining the size of the initial datum is maximised.

\begin{figure}[h!]
 \centering
 \includegraphics[width=0.35\textwidth,trim={1.8cm 0.8cm 1.25cm 0.5cm},clip]{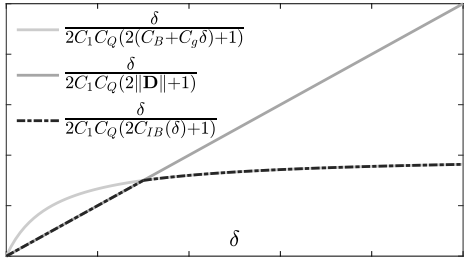}
 \caption{Bound on initial datum $\epsilon$ as a function of $\delta$}
 \label{fig:eps_vs_delta}
\end{figure}%

\begin{remark}
The ratio $\frac{\delta}{2C_1\eta}$, which is the part of $\varepsilon$ of interest here, is a monotonic increasing function for $\delta>0$. This can be seen in Fig.~\ref{fig:eps_vs_delta}, where this relationship has been portrayed for the case $\|\mathbf{D}\|>1$ (but equally valid otherwise). 
\end{remark}

\section{BOUND OPTIMISATION VIA SEMI-DEFINITE PROGRAMMING}
\label{sec:optimisation}
\subsection{Optimal Problem Definition}
As discussed above, the objective of the optimisation is to enlarge the region of attraction, which is bounded by the previously introduced constant $\epsilon$, by appropriate choice of the Lyapunov weighting and feedback matrices $Q(x)$ and $\kappa$. The proposed optimisation is independent of the beam's geometrical and material properties, and hence the constants $C_1$, $C_B$ and $C_g$ are assumed to be fixed known values. Therefore, the only constants which depend upon choice of $Q(x)$ are $C_Q$ and $C_\mathcal{S}$.
%
% \begin{figure}[h!]
%     \centering
%     \includegraphics[width=0.4\textwidth,trim={1.8cm 0.8cm 1.25cm 0.5cm},clip]{images/epsilon_delta.eps}
%     \caption{Bound on initial datum $\epsilon$ as a function of $\delta$}
%     \label{fig:eps_vs_delta}
% \end{figure}%
%
% By inspection of the dependence of $\epsilon$ on the previously defined constants, it is first observed that it is a monotonic increasing function for $\delta>0$. This is shown in Fig.~\ref{fig:eps_vs_delta}, where this relationship has been portrayed for the case $\|\mathbf{D}\|>1$ (but equally valid otherwise). Therefore, in order to maximise $\epsilon$, $\delta$ has to be picked as large as possible whilst respecting condition \eqref{eq:alpha}. 
Injecting the definition of $\delta$ provided by \eqref{eq:alpha} into \eqref{eq:analytical_bound}, and making use of the definitions \eqref{eq:eta}-\eqref{eq:def_CIB}, we obtain an upper bound for $\varepsilon$:
\begin{align}
\varepsilon \leq \frac{1}{4C_1 C_g \|L\|^2 \|L^{-1}\|} \left( \frac{C_\mathcal{S}}{C_Q^2} - \frac{ 2 \alpha}{C_Q^3} \right),
\label{eq:eps_bound}
\end{align}
where we recall that $\alpha>0$. Thus we will maximise the ratio $C_\mathcal{S}/{C_Q^2}$ in order to maximise the region of attraction estimate.
%Previous text:
%\textcolor{red}{If \eqref{eq:alpha} is taken as an equality, then an upper bound for $\epsilon$ can be defined, which would then scale as
%\begin{equation*}
 %   \epsilon\propto \frac{C_\mathcal{S}}{C_Q^2}.
%\end{equation*}}
It has been made obvious that for given structural and geometrical properties, the bound on the initial datum (i.e., the size of the region of attraction) is maximised by making $C_\mathcal{S}$ as large possible and choosing $C_Q$ as close to one as possible. These two choices are in direct competition, as dictated by inequalities \eqref{eq:MI}, and hence optimality is sought as proposed in the next section by employing a semi-definite programming approach.  

\subsection{Proposed Semi-Definite Program}
In the following, we consider the constant properties scenario, although the presented arguments can be directly generalised to the varying properties case if the property distributions are described, or can be well approximated, by polynomials of arbitrary degree. We rewrite the stability conditions \eqref{eq:MI} on the Lyapunov weighing matrix $Q(x)$ as conditions for non-negativity, suitable for semidefinite programming (SDP) tools, of the following matrix expressions, with $0<\epsilon_1,\,\epsilon_2\ll1$ (i.e., numerical artefacts used to obtain non-strict inequalities and can be made arbitrarily small):
\begin{subequations}
\label{eq:conditions_SDP}
\begin{align}
    Q-\epsilon_1\mathbf{I}_{12}&\succcurlyeq 0, \quad \forall x\in[0,\ell],\\
    -\tfrac{\mathrm{d}}{\mathrm{d}x}Q \mathbf{D} +QB + B^\top Q-\epsilon_2\mathbf{I}_{12}&\succcurlyeq 0, \quad \forall x\in[0,\ell],\\
    Q_{+}(\ell)-Q_{-}(\ell)&\succcurlyeq 0,\\
    Q_{-}(0)D-\kappa^\top Q_{+}(0)D\kappa &\succcurlyeq 0,\label{eq:cond_kappa}
\end{align}
\end{subequations}%
where the explicit dependence of the matrices $Q$ and $B$ on the spatial variable $x$ is not shown to alleviate notation.

%\textcolor{red}{We follow by defining each of the diagonal entries $q_i(x),\,i=1,\dots,12$, in $Q(x)=\mathrm{diag}(q_1(x),\dots,q_{12}(x))$ as polynomials in $x$ with arbitrary degree $n$, $q_i(x)=p_n(x)$.} 
We follow by defining each of the diagonal entries $\{q_i\}_{i=1}^{12}$ of $Q$ (see (5)) as polynomials in $x$ with arbitrary degree $n$, $q_i(x)=p_n(x)$. Since $\kappa$ is a decision variable, condition \eqref{eq:cond_kappa} introduces cubic terms (leading to a non-convex problem), however a Schur complement argument is employed to obtain an equivalent convex condition. Defining $\widetilde{\kappa}=\left(Q_{+}\left(0\right)D\right)^{\frac{1}{2}}{\kappa}$ and using the generalised \textit{s-procedure} \cite{BoydLMI} with the simple quadratic function $\mathcal{I}_\ell(x)=x(x-\ell)$, used to verify set containment, we finally write \eqref{eq:conditions_SDP} as
\begin{subequations}
\label{eq:feasibility_SDP}
\begin{align}
    Q-\epsilon_1\mathbf{I}_{12}+s_1(x)\mathcal{I}_\ell(x)\mathbf{I}_{12}&\succcurlyeq 0,\\
    -\tfrac{\mathrm{d}}{\mathrm{d}x}Q \mathbf{D} + QB + B^\top Q-\epsilon_2\mathbf{I}_{12}+s_2(x)\mathcal{I}_\ell(x)\mathbf{I}_{12}&\succcurlyeq 0,\\
    Q_{+}(\ell)-Q_{-}(\ell)&\succcurlyeq 0,\\
    \left[\begin{smallmatrix}
         \mathbf{I}_{6}& \widetilde{\kappa}  \\
         \widetilde{\kappa}^\top &  Q_{-}(0)D
    \end{smallmatrix}\right]&\succcurlyeq 0.
\end{align}
\end{subequations}%
Here, $s_1(x),\,s_2(x)\in\Sigma(x)$ are nonnegative polynomials ($\Sigma(x)$ is used to denote the set of all sum of squares polynomials in $x$). Hence, local exponential stability of \eqref{eq:IGEB} is guaranteed if \eqref{eq:feasibility_SDP} is feasible. 
%Feasibility of \eqref{eq:feasibility_SDP} has been rigorously demonstrated in \cite{RL2020} for particular weighting functions. In this paper, we rely on semi-definite programming to check feasibility within the space function of polynomials of arbitrary degree while seeking for optimality in the sense that has been introduced above.
Feasibility of \eqref{eq:feasibility_SDP} has been rigorously shown in \cite{RL2020} for particular weighting matrices constrained to be of the form $Q = \frac{1}{2}\mathrm{diag}(f \mathbf{M}, (2f(L) - f)\mathbf{M})$ for some $f\in C^1([0, \ell])$ fulfilling a series of properties.
The function $f$ found in \cite{RL2020} is not polynomial. One may also use \cite[(3.21)-(3.23)]{RL2020} to uncover an appropriate polynomial $f$ provided, however, that its degree is large enough (depending on the beam parameters and $\mathbf{E}$).
Here, we consider more general polynomial weighting matrices, whose diagonal entries may be chosen independently from one another, thus having the potential to fulfil the desired stabilisation task with smaller polynomial degree and to lead to a larger decay/region of attraction relying on semi-definite programming.

As previously discussed, maximising the region of attraction has been reduced to the maximisation of the ratio $C_\mathcal{S}/C_Q^2$. However, this defines a nonlinear objective function, which cannot be solved using standard semidefinite programming tools (only linear objective functions are supported). Instead, an iterative approach is considered, where at each iteration $C_\mathcal{S}$ is minimised for fixed maximum and minimum eigenvalues of $Q$ (which determine the constant $C_Q$). To achieve this, the following additions to \eqref{eq:feasibility_SDP} are considered to set up the SDP optimisation problem.

\emph{Maximising $C_\mathcal{S}$}. This can be achieved by imposing the following maximisation problem.
\begin{align*}
    &\max_\beta \beta\\
    \mathrm{s.t.}\,\, &-\tfrac{\mathrm{d}}{\mathrm{d}x}Q \mathbf{D}+ QB + B^\top Q-\beta \mathbf{I}_{12}\succcurlyeq0, \quad \forall x\in[0,\ell],\\
    &\beta\geq\epsilon_3.
\end{align*}

\emph{Bounds on $C_Q$}. This is enforced by adding the following constraints on the smallest and largest eigenvalues of $Q(x)$
% \begin{align*}
%     &Q(x)-\gamma I\succcurlyeq0, \quad \forall x\in[0,\ell],\\
%     &Q^{-1}(x)-\mu I\succcurlyeq0, \quad \forall x\in[0,\ell],\\
%     &0\leq\gamma<1,\\
%     &0\leq\mu<1.
% \end{align*}
% Or equivalently, and avoiding the appearance of the inverse $Q^{-1}$,
\begin{align*}
    &Q-\gamma \mathbf{I}_{12}\succcurlyeq0, \quad \forall x\in[0,\ell],\\
    &\nu \mathbf{I}_{12}-Q\succcurlyeq0, \quad \forall x\in[0,\ell],\\
    &0\leq\gamma\leq\nu.
\end{align*}

Gathering the previous additional maximisation sub-problem and constraints, the general SDP optimisation problem can now be defined

\begin{subequations}
\label{eq:optimal_prob_SDP}
\begin{align}
    &\min_{\beta,s_i(x)\in\Sigma(x),\,q_i(x),\,\widetilde{\kappa}} -\beta\\
    \mathrm{s.t.}\,\,\,\,&-\tfrac{\mathrm{d}}{\mathrm{d}x}Q \mathbf{D}\! +\! QB \!+\!B^\top Q\!-\beta \mathbf{I}_{12}\!+\!s_1(x)\mathcal{I}_\ell(x)\mathbf{I}_{12}\succcurlyeq 0,\\
    &Q_{+}(l)-Q_{-}(l)\succcurlyeq 0,\\
    &\left[\begin{smallmatrix}
         \mathbf{I}_{6}& \widetilde{\kappa}  \\
         \widetilde{\kappa}^\top&  Q_{-}(0)D
    \end{smallmatrix}\right]\succcurlyeq 0,\\
    &Q-\gamma \mathbf{I}_{12}+s_2(x)\mathcal{I}_\ell(x)\mathbf{I}_{12}\succcurlyeq0,\label{eq:max_eig_Q}\\
    &\nu\mathbf{I}_{12}-Q+s_{3}(x)\mathcal{I}_\ell(x)\mathbf{I}_{12}\succcurlyeq0,\label{eq:min_eig_Q}\\
    &\beta\geq\epsilon_1,\\
    &0\leq\gamma\leq\nu.\label{eq:evals_constraint}
\end{align}
\end{subequations}
Note that the positivity condition on the coefficients $Q(x)$ in \eqref{eq:conditions_SDP} has been made redundant by conditions \eqref{eq:max_eig_Q} and \eqref{eq:min_eig_Q}. 

Then, upon solution of \eqref{eq:optimal_prob_SDP} for fixed $\gamma,\,\nu$ a value of the ratio $C_\mathcal{S}/C_Q^2$ is obtained. Therefore, an iterative solution approach to find the optimal $C_\mathcal{S}/C_Q^2(\gamma,\nu)$ is employed, where each step will require solving for \eqref{eq:optimal_prob_SDP}.

% Note also that the optimisation problem below is feasible if and only if \eqref{eq:conditions_SDP} is feasible, for one can always pick suitable (although not optimal) $\beta,\,\gamma,\,\nu$ so that conditions of \eqref{eq:conditions_SDP} are recovered. Then, the outcome of \eqref{eq:optimal_prob_SDP} yields an optimised $\epsilon$ in some sense. Note that the solution of \eqref{eq:optimal_prob_SDP} is dependent on the chosen weights $a_\beta,\,a_\gamma,\,a_\nu$.

\section{NUMERICAL RESULTS}
\label{sec:results}

A numerical case to exemplify the proposed approach is performed on a beam with unitary structural and geometrical properties, that is, a beam with 
% $\mathbf{ M}=\mathbf{I}_6$ and $\mathbf{ C}=\mathbf{I}_6$, with appropriate units. 
the following characteristics.
\begin{table}[ht!]
    \centering
    \caption{Structural and geometrical properties}
    \label{tab:hale_discretisation}
    \begin{tabular}{lcl}
        \hline
        property & symbol & value \\
        \hline
        mass per unit length & $\rho a$ & $1\,\mathrm{kgm^{-1}}$\\
        area moments of inertia & $I_2,\,I_3$ & $1\,\mathrm{kgm}$\\
        axial stiffness & $Ea$ & $1\,\mathrm{N}$\\
        shear stiffness & $Ga$ & $1\,\mathrm{N}$\\
        torsional stiffness & $G(I_2+I_3)$ & $1\,\mathrm{Nm^2}$\\
        bending stiffness & $EI_2,\,EI_3$ & $1\,\mathrm{Nm^2}$\\
        polar moment of area correction & $k_1$ & 1\\
        shear correction factors & $k_2,\,k_3$ & 1\\
        beam length & $\ell$ & $1\,\mathrm{m}$\\ \hline
    \end{tabular}
\end{table}

The optimisation problem \eqref{eq:optimal_prob_SDP} has been implemented in Matlab using the toolbox SOSTOOLS~\cite{sostools} and the solver CDCS~\cite{Zheng}. Polynomials of degree $p=4$ are chosen to define the diagonal entries $q_i(x)$ of the weighing matrix $Q(x)$ and the non-negative functions $s_1(x),\,s_2(x),\,s_3(x)$ employed to verify set containment, which provide a good trade-off in terms of computational complexity of the underlying semi-definite problems.
\begin{figure}[t]
    \centering
    \includegraphics[width=0.38\textwidth,trim={0.25cm 0.95cm 1.25cm 0.35cm},clip]{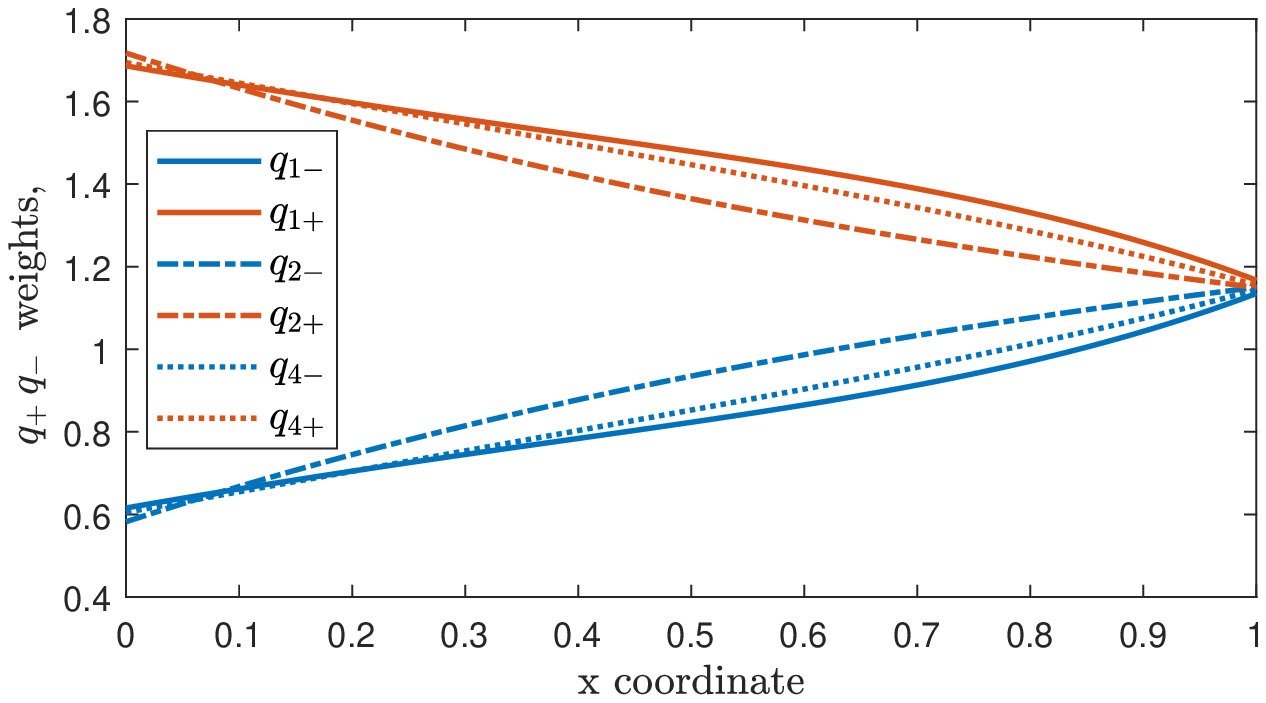}
    \includegraphics[width=0.38\textwidth,trim={0.25cm 0cm 1.25cm 0.25cm},clip]{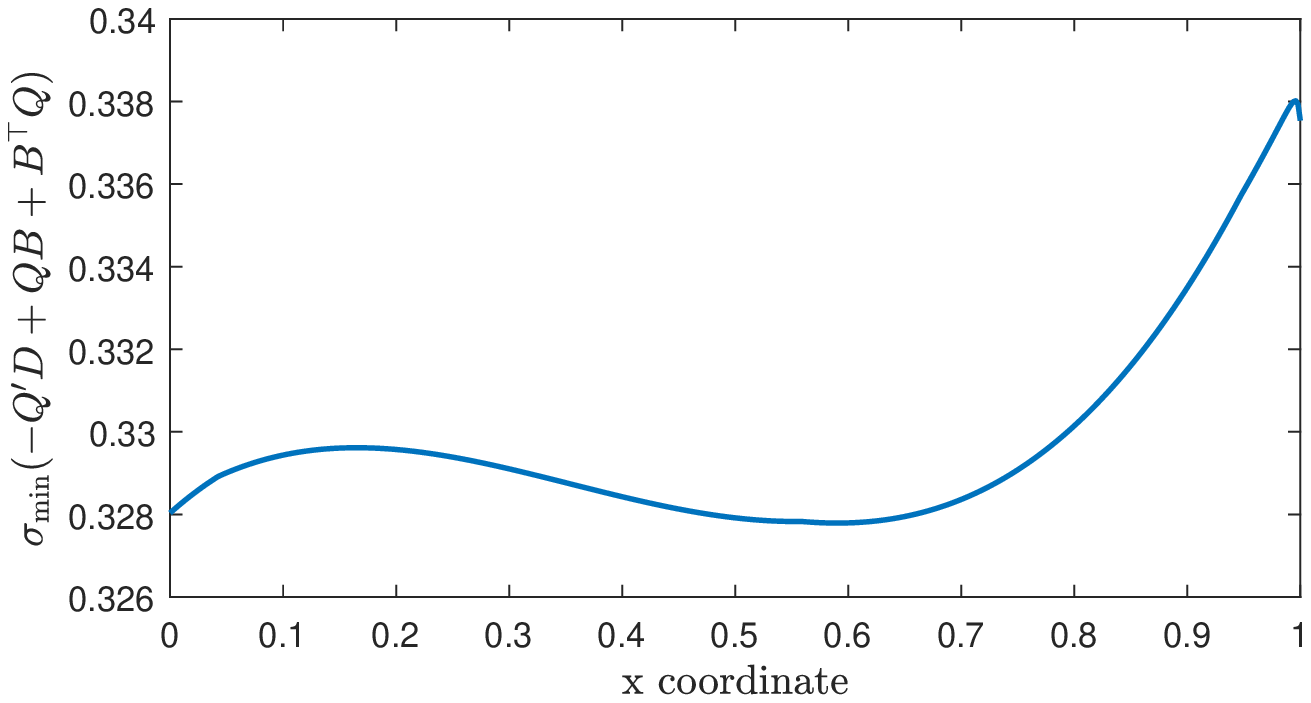}
     \caption{Diagonal entries of $Q(x)$ (top) and smallest eigenvalue of $\mathcal{S}$ (bottom) along spatial coordinate $x$ for the optimal pair $(\gamma,\nu)^\ast$}
    \label{fig:q_weights and sigma}
\end{figure}%
\begin{figure}[t]
    \centering
    \includegraphics[width=0.4\textwidth,trim={0.4cm 0.2cm 0cm 0.7cm},clip]{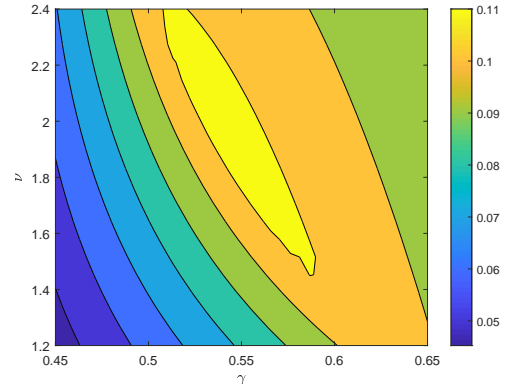}
     \caption{Contour plot of the ratio $C_\mathcal{S}/C_Q^2$ for varying $(\gamma,\nu)$}
    \label{fig:surf_contour}
\end{figure}%

The search for the optimal pair $(\gamma,\nu)^\ast$ which maximises the ratio $C_\mathcal{S}/C_Q^2$ is performed using Matlab's in-built optimisation function \emph{fminunc}. Despite the presence of constraint \eqref{eq:evals_constraint}, it has been found that, in practice, this condition is satisfied throughout the entire iterative process and hence the simpler unconstrained optimisation scenario is considered. The optimal pair has been found to be $(\gamma,\nu)^\ast=(0.5823,\,1.7173)$, resulting in constants $C_\mathcal{S}=0.3380$ and $C_Q=1.7173$ yielding a ratio $C_\mathcal{S}/C_Q^2=0.1146$. The entries of the Lyapunov functional weighing matrix \eqref{eq:Lyapunov_functional} along the spatial coordinate $x$ are displayed in Fig.~\ref{fig:q_weights and sigma}, where blue lines are used for the $q_{i_{-}}$ entries and red for the $q_{i_{+}}$. Coefficients $q_{i_{-}},\,q_{i_{+}}$ for $i=\{3,5,6\}$ are not included in the figure since they are found to be coincident to $q_{2_{-}}$ and $q_{2_{+}}$ to tolerance precision. This symmetry is attributed to the identical properties used in the two bending directions. A strong resemblance of the obtained shapes with the analytical functions obtained in \cite{RL2020} can be observed.

Fig.~\ref{fig:q_weights and sigma} also shows the smallest eigenvalue of the matrix $\mathcal{S}$ along the spatial coordinate $x$ for the optimal pair $(\gamma,\nu)^\ast$, whose minimum defines the constant $C_\mathcal{S}$. The obtained optimal feedback matrix is
% $\kappa$ (and $\overline{\kappa}$) are
% \begin{equation*}
% \kappa=10^{-2}\cdot\textrm{diag}(-0.39,\,6.81,\,7.82,\,-0.89,\,7.93,\,7.95),
% \end{equation*}
\begin{equation*}
\bar{\kappa}=\textrm{diag}(1.0079,\,0.8725,\,0.8549,\,1.4397,\,0.8530,\,0.8528).
\end{equation*}

%\textcolor{red}{Once the constants $C_\mathcal{S}$ and $C_Q$ are obtained from the optimisation solution, the value of the constant $\delta$ is computed by taking \eqref{eq:alpha} as an equality.} 
% Once the constants $C_\mathcal{S}$ and $C_Q$ are obtained from the optimisation solution,
The value of the constant $\delta$ is given by \eqref{eq:alpha} together with our choice for $0<\alpha<\smash{\frac{1}{2}}C_\mathcal{S}C_Q$. Given that the structural properties are known and fixed, the constant $\eta$ is readily available from \eqref{eq:eta} and hence the bound on the region of attraction is finally obtained from \eqref{eq:eps_bound}. In the limit $\alpha=0$, our numerically found bound is $\epsilon<0.0095$. The sensitivity of the numerical results to the polynomial degree has been observed to be rather low, since the optimal pair $(\gamma,\nu)^\ast$ produce very similar solutions to \eqref{eq:optimal_prob_SDP} for both lower and higher polynomial degrees (2 and 6).

The ratio $C_\mathcal{S}/C_Q^2$, and consequently the bound on the region of attraction $\epsilon$, shows a strong dependence on the choice of $Q(x)$ and the imposed constraints on $Q(x)$ introduced by the constants $(\gamma,\nu)$, which justifies the need for exploration of a suitable weighing matrix to enlarge the region from which solutions are expected to decay exponentially. This is clearly shown in Fig.~\ref{fig:surf_contour}, where the ratio $C_\mathcal{S}/C_Q^2$ has been plotted for a range of $(\gamma,\,\nu)$, obtained through a sweep over 20 different values for each constant around the optimum.

\section{CONCLUSIONS}
An optimisation method to obtain sharper bounds on the region of attraction for the exponential stabilisation of geometrically exact beams via boundary feedback has been demonstrated. This strategy, based on semi-definite programming, introduces a systematic approach to select Lyapunov functionals which gives a better insight on the size of initial datum leading to exponential decay of solutions. The proposed method can equally be applied to more general, spatially varying mass and flexibility matrices, relying on an equivalent theoretical proof, with the (physical) constraint that they are symmetric positive definite matrices.

This methodology is, however, still restricted to the small initial datum scenario, since the established stability results rely on the dominance of the linear subsystem over the nonlinear terms. Besides, control laws necessary to achieve stability under these results are required to provide feedback on all degrees of freedom. A future line of investigation to overcome these limitations and to explore for more general stability or boundedness results is to consider directly the inequality $\partial_t\mathcal{L}+\alpha\mathcal{L}\leq0$, for some $\alpha>0$. The use of semi-definite programming tools on this full (nonlinear) expression offers a viable alternative to comprehend and estimate the role of nonlinear couplings in stability analysis, which is an intractable task if only analytical tools are considered.

\addtolength{\textheight}{-12cm}   % This command serves to balance the column lengths

%%%%%%%%%%%%%%%%%%%%%%%%%%%%%%%%%%%%%%%%%%%%%%%%%%%%%%%%%%%%%%%%%%%%%%%%%%%%%%%%

%%%%%%%%%%%%%%%%%%%%%%%%%%%%%%%%%%%%%%%%%%%%%%%%%%%%%%%%%%%%%%%%%%%%%%%%%%%%%%%%


% \section*{ACKNOWLEDGMENT}

%%%%%%%%%%%%%%%%%%%%%%%%%%%%%%%%%%%%%%%%%%%%%%%%%%%%%%%%%%%%%%%%%%%%%%%%%%%%%%%%

\bibliographystyle{IEEEtran} 
\bibliography{IEEEabrv.bib,mybib}

\end{document}